\documentstyle[12pt,amsfonts,twoside,leqno]{amsart}

\newtheorem{Theorem}{Theorem}[section]
\newtheorem{Proposition}[Theorem]{Proposition}
\newtheorem{Lemma}[Theorem]{Lemma}
\newtheorem{Corollary}[Theorem]{Corollary}

\theoremstyle{remark}
\newtheorem{Remark}[Theorem]{Remark}

\theoremstyle{definition}

\begin{document}
\bibliographystyle{plain}
\title{Conductor inequalities and criteria for Sobolev-Lorentz two-weight inequalities}
\author{\c Serban Costea \and Vladimir Maz'ya}

\address{\c Serban Costea: secostea@@math.mcmaster.ca, McMaster University,
Department of Mathematics and Statistics, 1280 Main Street West, Hamilton, Ontario L8S 4K1, Canada \\ \and
Fields Institute for Research in Mathematical Sciences, 222 College Street, Toronto, Ontario M5T 3J1, Canada}

\address{Vladimir Maz'ya: vlmaz@@math.ohio-state.edu, vlmaz@@liv.ac.uk, vlmaz@@mai.liu.se, Department of Mathematics, The Ohio State University, 231 W 18th Ave, Columbus, OH 43210, USA \\ Department of Mathematical Sciences, M O Building,  University of Liverpool, Liverpool L69 3BX, UK \\ \and
Department of Mathematics, Link\"{o}ping University, SE-581 83 Link\"{o}ping, Sweden}

\maketitle
\centerline
{\it In memory of S. L. Sobolev}

\bigskip

\begin{abstract}

  In this paper we present  integral conductor inequalities connecting the Lorentz $p,q$-(quasi)norm of a gradient of a function to a one-dimensional integral of the $p,q$-capacitance of the conductor between two level surfaces of the same function. These inequalities
  generalize an  inequality obtained by the second author in the case of the Sobolev norm. Such conductor inequalities lead to necessary and sufficient conditions for Sobolev-Lorentz type inequalities involving two arbitrary measures.

\end{abstract}

\maketitle
\section{\bf {Introduction}}

  During the last decades Sobolev-Lorentz function spaces, which include classical Sobolev spaces, attracted attention not only as an interesting
  mathematical object, but also as a tool for a finer tuning of properties of solutions to partial differential equations.
  (See \cite{Alb}, \cite{AFT1}, \cite{AFT2}, \cite{BBGGPV}, \cite{Cia}, \cite{CP}, \cite{Cos}, \cite{DHM}, \cite{HL}, \cite{KKM}, \cite{ST}, {\it{et al}}.)

  In the present paper we generalize the inequality
 \begin{equation} \label{stronger p Sobolev cap}
 \int_{0}^{\infty} {\mathrm{cap}}_p(\overline{M_{at}}, M_t) d(t^p) \le c(a,p) \int_{\Omega} |\nabla f|^p \, dx
 \end{equation}
 to Sobolev-Lorentz spaces. Here $f\in Lip_{0}(\Omega)$, i.e. $f$ is an arbitrary Lipschitz function compactly supported in the open set $\Omega \subset {\mathbf{R}}^n,$ while
  $M_t$ is the set $\{ x \in \Omega: |f(x)|>t \}$ with $t > 0.$
 Inequality (\ref{stronger p Sobolev cap}) was obtained in \cite{M1}. (See also \cite[Chapter 2]{M3}.) It has various extensions and applications
 to the theory of Sobolev-type spaces on domains in ${\mathbf{R}}^n,$ Riemannian manifolds, metric and topological spaces, to linear and nonlinear partial differential equations, Dirichlet forms, and Markov processes etc. (See \cite{Ad}, \cite{AH}, \cite{AP}, \cite{AX1}, \cite{AX2}, \cite{Ai},
 \cite{CS}, \cite{DKX}, \cite{Dah}, \cite{Fi}, \cite{FU1}, \cite{FU2}, \cite{Gr}, \cite{Haj}, \cite{Han}, \cite{HMV}, \cite{Ka}, \cite{Ko1}, \cite{Ko2}, \cite{Mal}, \cite{M1}, \cite{M2}, \cite{M4}, \cite{M5}, \cite{MN}, \cite{MP}, \cite{Ne}, \cite{Ra}, \cite{Ta},  \cite{V1}, \cite{V2}, \cite{Vo}, {\it et al}).
 In the sequel, we prove the inequalities
 \begin{equation} \label{stronger cap Sobolev-Lorentz q le p}
 \int_{0}^{\infty} {\mathrm{cap}}(\overline{M_{at}}, M_t) d(t^{p}) \le c(a,p,q) || \nabla f ||_{L^{p,q}(\Omega, m_n; {\mathbf{R}}^n)}^{p} \mbox{ when $1 \le q \le p$}
 \end{equation}
 and
 \begin{equation} \label{stronger cap Sobolev-Lorentz p le q}
 \int_{0}^{\infty} {\mathrm{cap}}_{p,q}(\overline{M_{at}}, M_t)^{q/p} d(t^{q}) \le c(a,p,q) || \nabla f ||_{L^{p,q}(\Omega, m_n; {\mathbf{R}}^n)}^{q} \mbox{ when $p < q <\infty$}
 \end{equation}
for all $f\in Lip_{0}(\Omega)$.

 The proof of (\ref{stronger cap Sobolev-Lorentz q le p}) and (\ref{stronger cap Sobolev-Lorentz p le q}) is based on the superadditivity of
 the $p,q$-capacitance, also justified in this paper.

 From (\ref{stronger cap Sobolev-Lorentz q le p}) and (\ref{stronger cap Sobolev-Lorentz p le q}) we derive necessary and sufficient conditions for certain two-weight inequalities involving Sobolev-Lorentz norms, generalizing results obtained in \cite{M4} and \cite{M5}. Specifically, let $\mu$ and $\nu$ be two locally finite nonnegative measures on $\Omega$ and let $p,q,r,s$ be real numbers such that $1<s \le \max (p,q) \le r<\infty$
 and $q \ge 1.$ We characterize the inequality
 \begin{equation} \label{embedding of Sobolev-Lorentz type space into Lebesgue spaces}
 ||f||_{L^{r,\max(p,q)}(\Omega, \mu)} \le A \left( || \nabla f ||_{L^{p,q}(\Omega, m_n; {{\mathbf{R}}^n})} + ||f||_{L^{s, \max(p,q)}(\Omega, \nu)} \right)
 \end{equation}
 restricted to functions $f \in Lip_{0}(\Omega)$ by requiring the condition
 \begin{equation} \label{two measure inequality for general n}
 \mu(g)^{1/r} \le K ({\mathrm{cap}}_{p,q}(\overline{g}, G)^{1/p} + \nu(G)^{1/s})
 \end{equation}
 to be valid for all open bounded sets $g$ and $G$ subject to $\overline{g} \subset G,$ $\overline{G} \subset \Omega.$
 When $n=1$ inequality (\ref{embedding of Sobolev-Lorentz type space into Lebesgue spaces}) becomes
 \begin{equation} \label{embedding of Sobolev-Lorentz type space into Lebesgue spaces n is 1}
 ||f||_{L^{r,\max(p,q)}(\Omega, \mu)} \le A \left( ||f'||_{L^{p,q}(\Omega, m_1)} + ||f||_{L^{s, \max(p,q)}(\Omega, \nu)} \right).
 \end{equation}
 The requirement that (\ref{embedding of Sobolev-Lorentz type space into Lebesgue spaces n is 1}) be valid for all functions
 $f \in Lip_{0}(\Omega)$ when $n=1$ is shown to be equivalent to the condition
 \begin{equation} \label{two measure inequality for centered intervals when n equal 1}
 \mu(\sigma_d(x))^{1/r} \le K (\tau^{(1-p)/p} + \nu(\sigma_{d+\tau}(x))^{1/s})
 \end{equation}
 whenever $x,$ $d$ and $\tau$ are such that $\overline{\sigma_{d+\tau}(x)} \subset \Omega.$
 Here and throughout the paper $\sigma_{d}(x)$ denotes the open interval $(x-d, x+d)$ for every $d>0.$

\section{\bf{Preliminaries}}

 Let us introduce some notation, to be used in the sequel. By $\Omega$ we denote a nonempty open subset of ${\mathbf{R}}^n,$ whereas $m_n$
 stands for the Lebesgue $n$-measure in ${{\mathbf{R}}}^n,$ where $n \ge 1$ is integer.
 For a Lebesgue measurable $u: \Omega \rightarrow {\mathbf{R}},$ $\mbox {supp } u$ is the
 smallest closed set such that $u$ vanishes outside $\mbox {supp } u.$
 We also define
 \begin{eqnarray*}
 Lip(\Omega)&=&\{ \varphi: \Omega \rightarrow {\mathbf{R}}: \varphi \mbox{ is Lipschitz} \}\\
 Lip_{0}(\Omega)&=&\{ \varphi: \Omega \rightarrow {\mathbf{R}}: \varphi \mbox{ is Lipschitz and with compact support in $\Omega$} \}.
 \end{eqnarray*}
 If $\varphi \in Lip(\Omega),$ we write $\nabla \varphi$ for the gradient of $\varphi.$ This notation makes sense, since
 by Rademacher's theorem (\cite[Theorem 3.1.6]{Fed}) every Lipschitz function
 on $\Omega$ is $m_n$-a.e. differentiable.

 Throughout this section we will assume that $m \ge 1$ is a positive integer and that $(\Omega, \mu)$
 is a measure space.
 Let $f:\Omega \rightarrow {\mathbf{R}}^n$ be a $\mu$-measurable function.
 We define $\mu_{[f]},$ the \textit{distribution function}
 of $f$ as follows (see \cite[Definition II.1.1]{BS}):
 $$\mu_{[f]}(t)=\mu(\{x \in \Omega: |f(x)| > t \}), \qquad t \ge 0.$$
 We define $f^{*},$  the \textit{nonincreasing rearrangement} of $f$ by
$$f^{*}(t)=\inf\{v: \mu_{[f]}(v) \le t \}, \quad t \ge 0.$$
(See \cite[Definition II.1.5]{BS}.)
 We notice that $f$ and $f^{*}$ have the same distribution function.
 Moreover, for every positive $\alpha$ we have $$(|f|^{\alpha})^{*}=(|f|^{*})^{\alpha}$$
 and if $|g|\le |f|$ a.e. on $\Omega,$ then $g^{*}\le f^{*}.$
 (See \cite[Proposition II.1.7]{BS}.)
 We also define $f^{**}$, the \textit{maximal function} of $f^{*}$ by
  $$f^{**}(t)=m_{f^{*}}(t)=\frac{1}{t} \int_{0}^{t} f^{*}(s) ds, \quad t >0.$$
 (See \cite[Definition II.3.1]{BS}.)

Throughout this paper, we will denote  the H\"{o}lder
conjugate of $p \in [1,\infty]$ by $p'$.

\smallskip

The \textit{Lorentz space} $L^{p,q}(\Omega, \mu; {\mathbf{R}}^n),$
$1<p<\infty,$ $1\le q\le \infty,$ is defined as follows:
$$L^{p,q}(\Omega, \mu; {\mathbf{R}}^n)= \{f: \Omega \rightarrow {\mathbf{R}}^n: f \mbox { is $\mu$-measurable and }
||f||_{L^{p,q}(\Omega, \mu; {\mathbf{R}}^n)}<\infty\},$$
where
$$||f||_{L^{p,q}(\Omega,\mu; {\mathbf{R}}^n)}=||\,|f|\,||_{p,q}=\left\{ \begin{array}{lc}
\left( \displaystyle{\int_{0}^{\infty} (t^{1/p}f^{*}(t))^q \, \frac{dt}{t}}
\right)^{1/q} & 1 \le q < \infty \\
\\
\sup\limits_{t>0} t \mu_{[f]}(t)^{1/p}=\sup\limits_{s>0}
s^{1/p} f^{*}(s) & q=\infty.
\end{array}
\right.
$$
(See \cite[Definition IV.4.1]{BS} and \cite[p.\ 191]{SW}.) We omit ${\mathbf{R}}^n$ in the notation of function spaces for the scalar case, i.e. for $n=1.$

\smallskip

If $1 \le q\le p,$ then $||\cdot||_{L^{p,q}(\Omega, \mu; {\mathbf{R}}^n)}$
represents a norm, but for $p < q \le \infty$ it represents a
quasinorm, equivalent to the norm $||\cdot||_{L^{(p,q)}(\Omega, \mu;
{\mathbf{R}}^n)},$ where
$$||f||_{L^{(p,q)}(\Omega, \mu; {\mathbf{R}}^n)}=||\,|f|\,||_{(p,q)}=\left\{ \begin{array}{lc}
\left(\displaystyle{ \int_{0}^{\infty} (t^{1/p}f^{**}(t))^q \, \frac{dt}{t} }\right)^{1/q} & 1 \le q < \infty \\
\\
\sup\limits_{t>0} t^{1/p} f^{**}(t) & q=\infty.
\end{array}
\right.
$$
(See \cite[Definition IV.4.4]{BS}.) Namely, from \cite[Lemma
IV.4.5]{BS} we have that
$$||\,|f|\,||_{L^{p,q}(\Omega,\mu)} \le ||\,|f|\,||_{L^{(p,q)}(\Omega,\mu)} \le p' ||\,|f|\,||_{L^{p,q}(\Omega,\mu)}$$
for every $q \in [1, \infty]$ and every $\mu$-measurable function
$f:\Omega \rightarrow {\mathbf{R}}^n.$

It is known that $(L^{p,q}(\Omega,\mu; {\mathbf{R}}^n),
||\cdot||_{L^{p,q}(\Omega,\mu; {\mathbf{R}}^n)})$ is a Banach space for
$1\le q \le p,$ while $(L^{p,q}(\Omega,\mu; {\mathbf{R}}^n),
||\cdot||_{L^{(p,q)}(\Omega,\mu; {\mathbf{R}}^n)})$ is a Banach space
for $1<p< \infty,$ $1\le q \le \infty.$

\begin{Remark} \label{relation between Lpr and Lps} It is also known (see
\cite[Proposition IV.4.2]{BS}) that for every $p \in (1,\infty)$ and
$1\le r<s\le \infty$ there exists a constant $C(p,r,s)$ such that
\begin{equation}\label{relation between Lpr and Lps norm}
||\,|f|\,||_{L^{p,s}(\Omega, \mu)} \le C(p,r,s)
||\,|f|\,||_{L^{p,r}(\Omega, \mu)}
\end{equation}
for all measurable functions $f \in L^{p,r}(\Omega, \mu;{\mathbf{R}}^n)$
and all integers $n \ge 1.$ In particular, the embedding
$L^{p,r}(\Omega, \mu;{\mathbf{R}}^n) \hookrightarrow
L^{p,s}(\Omega, \mu;{\mathbf{R}}^n)$ holds.
\end{Remark}

\subsection{\bf{The subadditivity and superadditivity of the Lorentz quasinorms}}

In the second part of this paper, we will prove a few results by relying on the superadditivity
of the Lorentz $p,q$-quasinorm. Therefore we recall the known results and present new results concerning the superadditivity
and the subadditivity of the Lorentz $p,q$-quasinorm.

The superadditivity of the Lorentz $p,q$-norm in the case $1 \le q \le p$ was stated in \cite[Lemma 2.5]{CHK}.

\begin{Proposition} \label{superadd q le p} (See \cite[Lemma 2.5]{CHK}.)
Let $(\Omega, \mu)$ be a measure space. Suppose $1 \le q \le p.$ Let $\{E_i\}_{i \ge 1}$ be
a collection of pairwise disjoint measurable subsets of $\Omega$ with $E_{0}=\cup_{i \ge 1} E_i$ and let $f \in L^{p,q}(\Omega, \mu).$ Then
\begin{equation*}
 \sum_{i \ge 1} || \chi_{E_i} f||_{L^{p,q}(\Omega, \mu)}^p \le || \chi_{E_0} f||_{L^{p,q}(\Omega, \mu)}^p.
\end{equation*}
\end{Proposition}

We obtain a similar result concerning the superadditivity in the case $1<p < q<\infty.$

\begin{Proposition} \label{superadd p le q} Let $(\Omega, \mu)$ be a measure space. Suppose $1<p < q<\infty.$ Let $\{E_i\}_{i \ge 1}$ be
a collection of pairwise disjoint measurable subsets of $\Omega$ with $E_{0}=\cup_{i \ge 1} E_i$ and let $f \in L^{p,q}(\Omega, \mu).$ Then
\begin{equation*}
 \sum_{i \ge 1} || \chi_{E_i} f||_{L^{p,q}(\Omega, \mu)}^q \le || \chi_{E_0} f||_{L^{p,q}(\Omega, \mu)}^q.
\end{equation*}
\end{Proposition}

\begin{pf} For every $i=0,1,2,\ldots$ we let $f_{i}= \chi_{E_{i}} f,$ where $\chi_{E_i}$ is the characteristic
function of $E_{i}.$ We can assume without loss of generality that all the functions $f_i$ are nonnegative.
We have (see \cite[Proposition 2.1]{KKM})
$$||f_i||_{L^{p,q}(\Omega, \mu)}^{q}=p \int_{0}^{\infty} s^{q-1}
{\mu_{[f_i]}(s)}^{q/p} ds,$$ where $\mu_{[f_i]}$ is the distribution function of $f_i, i=0,1,2,\ldots.$
From the definition of $f_0$ we have
\begin{equation} \label{the distrib fn of f0 equals the sum of distrib fn fi}
\mu_{[f_0]}(s)=\sum_{i \ge 1} \mu_{[f_i]}(s) \mbox{ for every $s>0,$}
\end{equation}
which implies, since
$1 < p < q<\infty,$ that
$$\mu_{[f_0]}(s)^{q/p} \ge \sum_{i \ge 1} \mu_{[f_i]}(s)^{q/p} \mbox{ for every $s>0.$}$$ This yields
\begin{eqnarray*}
||f_0||_{L^{p,q}(\Omega, \mu)}^{q}&=&p \int_{0}^{\infty} s^{q-1}
{\mu_{[f_0]}(s)}^{q/p} ds \ge p \int_{0}^{\infty} s^{q-1}
(\sum_{i \ge 1} {\mu_{[f_i]}(s)}^{q/p} ) ds\\
 &=& \sum_{i \ge 1} p \int_{0}^{\infty} s^{q-1}
 {\mu_{[f_i]}(s)}^{q/p}  ds = \sum_{i \ge 1} ||f_i||_{L^{p,q}(\Omega, \mu)}^{q}.
\end{eqnarray*}
This finishes the proof of the superadditivity in the case $1< p < q<\infty.$

\end{pf}

We have a similar result for the subadditivity of the Lorentz $p,q$-quasinorm.
When $1 < p < q \le \infty$ we obtain a result that generalizes \cite[Theorem 2.5]{Cos}.

\begin{Proposition} \label{subadd p le q} Let $(\Omega, \mu)$ be a measure space. Suppose $1 < p < q \le \infty.$ Let $\{E_i\}_{i \ge 1}$ be
a collection of pairwise disjoint measurable subsets of $\Omega$ with $E_{0}=\cup_{i \ge 1} E_i$ and let
$f \in L^{p,q}(\Omega, \mu).$ Then

 $$\sum_{i \ge 1} || \chi_{E_i} f||_{L^{p,q}(\Omega, \mu)}^p \ge || \chi_{E_0} f||_{L^{p,q}(\Omega, \mu)}^p.$$

\end{Proposition}

\begin{pf} Without loss of generality we can assume that all the functions $f_i=\chi_{E_i} f$
are nonnegative. We have to consider two cases, depending
on whether $p<q<\infty$ or $q=\infty.$

Suppose $p<q<\infty.$ We have (see \cite[Proposition 2.1]{KKM})
$$||f_i||_{L^{p,q}(\Omega, \mu)}^{p}=\left(p \int_{0}^{\infty} s^{q-1}
{\mu_{[f_i]}(s)}^{q/p} ds\right)^{p/q},$$ where
$\mu_{[f_i]}$ is the distribution function of $f_i$ for
$i=0,1,2, \ldots.$ From (\ref{the distrib fn of f0 equals the sum of distrib fn fi}) we obtain
\begin{eqnarray*}
||f_0||_{L^{p,q}(\Omega, \mu)}^{p} & = & \left(p \int_{0}^{\infty}
s^{q-1}  \mu_{[f_0]}(s)^{q/p} ds\right)^{p/q} \le \sum_{i \ge 1} \left(p \int_{0}^{\infty} s^{q-1}
{\mu_{[f_i]}(s)}^{q/p} ds\right)^{p/q} \\
& = &\sum_{i \ge 1} ||f_i||_{L^{p,q}(\Omega, \mu)}^{p}.
\end{eqnarray*}

Suppose now $q=\infty.$ From
(\ref{the distrib fn of f0 equals the sum of distrib fn fi}) we obtain
\begin{equation*}
s^p \, \mu_{[f_0]}(s)  = \sum_{i \ge 1} (s^p \, \mu_{[f_i]}(s)) \mbox{ for every $s>0,$}
\end{equation*}
which implies
\begin{equation}\label{p, inf equality 1}
s^p \, \mu_{[f_0]}(s) \le \sum_{i \ge 1} ||f_i||_{L^{p, \infty}(\Omega, \mu)}^p \mbox{ for every $s>0.$}
\end{equation}
By taking the supremum over all $s>0$ in
(\ref{p, inf equality 1}), we get the desired conclusion. This finishes the proof.

\end{pf}

\section{\bf{Sobolev-Lorentz $p,q$-capacitance}}

Suppose $1<p<\infty$ and $1 \le q \le \infty.$ Let $\Omega \subset {\mathbf{R}}^n$ be
an open set, $n \ge 1.$ Let $K \subset \Omega$ be compact. The Sobolev-Lorentz
$p,q$-capacitance of the conductor $(K, \Omega)$ is denoted by
$${\mathrm{cap}}_{p,q}(K, \Omega)=\inf \, \{||\nabla u||_{L^{p,q}(\Omega, m_n; {\mathbf{R}}^n)}^{p}:
u \in W(K, \Omega)\},$$ where
\begin{equation*}
W(K, \Omega)=\{ u \in Lip_{0}(\Omega): u \ge 1\ \mbox{ in a
neighborhood of } K\}.
\end{equation*}
We call $W(K,\Omega)$ the \textit{set of admissible functions for
the conductor $(K, \Omega).$}

Since $W(K, \Omega)$ is closed under truncations from below by $0$ and from above by
$1$ and since these truncations do not increase the $p,q$-quasinorm whenever $1<p<\infty$ and
$1 \le q\le \infty$, it follows that we can choose only functions $u \in W(K,
\Omega)$ that satisfy $0 \le u \le 1$ when computing the $p,q$-capacitance of the conductor $(K, \Omega).$

\begin{Lemma} \label{same capacity if Lip and zero on bdry considered}
If $\Omega$ is bounded, then we get the same $p,q$-capacitance for the
conductor $(K, \Omega)$ if we restrict ourselves to a bigger set, namely
\begin{equation*}
 W_{1}(K,\Omega)=\{ u \in Lip(\Omega) \cap C(\overline{\Omega}) :  u \ge 1
 \mbox{ on } K \mbox{ and } u=0 \mbox{ on $\partial \Omega$} \}.
 \end{equation*}
\end{Lemma}

\begin{pf} Let $u \in W_{1}(K, \Omega).$ We can assume without loss of generality that $0 \le u \le 1.$ Moreover, we can also assume that $u=1$ in an open neighborhood $U$ of $K.$ Let $\widetilde{U}$ be an open neighborhood of $K$ such that
$\widetilde{U} \subset \subset U.$ We choose a cutoff Lipschitz function $\eta,$ $0 \le \eta \le 1$ such that $\eta=1$ on
$\Omega \setminus U$ and $\eta=0$ on $\widetilde{U}.$ We notice that $1-\eta(1-u)=u.$ We also notice that there exists a sequence of functions $\varphi_j \in Lip_{0}(\Omega)$ such that
$$\lim_{j \rightarrow \infty} (||\varphi_j-u||_{L^{p+1}(\Omega, m_n)} + ||\nabla \varphi_j-\nabla u||_{L^{p+1}(\Omega, m_n; {\mathbf{R}}^n)})= 0.$$  Without loss of generality the sequence $\varphi_j$ can be chosen such that $\varphi_j \rightarrow u$ and $\nabla \varphi_j \rightarrow \nabla u$ pointwise a.e.\ in $\Omega.$ Then $\psi_{j}=1-\eta(1-\varphi_{j})$ is a sequence
 belonging to $W(K,\Omega)$ and
 $$\lim_{j \rightarrow \infty} (||\psi_j-u||_{L^{p+1}(\Omega, m_n)} + ||\nabla \psi_j-\nabla u||_{L^{p+1}(\Omega, m_n; {\mathbf{R}}^n)})=0.$$
 This, H\"{o}lder's inequality for Lorentz spaces, and the behaviour of the Lorentz $p,q$-quasinorm in $q$ yield
 $$\lim_{j \rightarrow \infty} (||\psi_j-u||_{L^{p,q}(\Omega, m_n)} + ||\nabla \psi_j-\nabla u||_{L^{p,q}(\Omega, m_n; {\mathbf{R}}^n)})=0.$$
 The desired conclusion follows.

\end{pf}

\subsection{\bf{Basic properties of the $p,q$-capacitance}}

Usually, a capacitance is a monotone and subadditive set function. The
following theorem will show, among other things, that this is true
in the case of the $p,q$-capacitance. We follow \cite{Cos} for (i)-(vi).
In addition we will prove some superadditivity properties of the $p,q$-capacitance.

\begin{Theorem}\label{Cap Thm 1<q leq infty} Suppose $1<p<\infty$ and $1 \le q \le \infty.$
Let $\Omega \subset {\mathbf{R}}^n$ be open. The set function $K
\mapsto {\mathrm{cap}}_{p,q}(K, \Omega),$ $K \subset \Omega,$ $K$
compact, enjoys the following properties:

\par {\rm{(i)}}  If $K_{1} \subset K_{2},$ then ${\mathrm{cap}}_{p,q}(K_{1},
\Omega) \le {\mathrm{cap}}_{p,q}(K_{2}, \Omega).$

\par {\rm{(ii)}} If $\Omega_{1} \subset \Omega_{2}$ are open and $K$ is a compact
subset of $\Omega_{1},$ then $${\mathrm{cap}}_{p,q}(K, \Omega_{2})
\le {\mathrm{cap}}_{p,q}(K, \Omega_{1}).$$

\par {\rm{(iii)}} If $K_{i}$ is a decreasing sequence of compact subsets of
$\Omega$ with $K=\bigcap_{i=1}^{\infty} K_{i},$ then
$${\mathrm{cap}}_{p,q}(K, \Omega)=\lim_{i \rightarrow \infty}
{\mathrm{cap}}_{p,q}(K_{i}, \Omega).$$

\par {\rm{(iv)}} If $\Omega_i$ is an increasing sequence of open
sets with $\bigcup_{i=1}^{\infty} \Omega_{i}=\Omega$ and $K$ is a
compact subset of $\Omega_1,$ then
$${\mathrm{cap}}_{p,q}(K, \Omega)=\lim_{i \rightarrow \infty}
{\mathrm{cap}}_{p,q}(K, \Omega_i).$$

\par {\rm{(v)}} Suppose $p \le q \le \infty.$ If $K=\bigcup_{i=1}^{k} K_{i} \subset \Omega$ then
$${\rm{cap}}_{p,q}(K, \Omega) \le \sum_{i=1}^{k} {\rm{cap}}_{p,q}(K_{i}, \Omega),$$
where $k \ge 1$ is a positive integer.

\par {\rm{(vi)}} Suppose $1 \le q < p.$ If $K=\bigcup_{i=1}^{k} K_{i} \subset \Omega$ then
$${\rm{cap}}_{p,q}(K, \Omega)^{q/p} \le \sum_{i=1}^{k} {\rm{cap}}_{p,q}(K_{i}, \Omega)^{q/p},$$
where $k \ge 1$ is a positive integer.

\par {{\rm(vii)}} Suppose $1 \le q \le p.$  Suppose $\Omega_i, \ldots, \Omega_k$ are $k$ pairwise disjoint open sets and $K_i$
are compact subsets of $\Omega_i$ for $i=1, \ldots, k.$ Then
 $${\mathrm{cap}}_{p,q}(\cup_{i=1}^k K_i, \cup_{i=1}^k \Omega_i) \ge \sum_{i=1}^{k} {\mathrm{cap}}_{p,q}(K_i, \Omega_i).$$

\par {{\rm(viii)}} Suppose $p < q < \infty.$  Suppose $\Omega_i, \ldots, \Omega_k$ are $k$ pairwise disjoint open sets and $K_i$
are compact subsets of $\Omega_i$ for $i=1, \ldots, k.$ Then
 $${\mathrm{cap}}_{p,q}(\cup_{i=1}^k K_i, \cup_{i=1}^k \Omega_i)^{q/p} \ge \sum_{i=1}^{k} {\mathrm{cap}}_{p,q}(K_i, \Omega_i)^{q/p}.$$

\par {\rm{(ix)}} Suppose $1 \le q < \infty.$ If $\Omega_1$ and $\Omega_2$ are two disjoint open sets and $K \subset \Omega_1,$
then $${\mathrm{cap}}_{p,q}(K, \Omega_1 \cup \Omega_2)={\mathrm{cap}}_{p,q}(K, \Omega_1).$$

\end{Theorem}

\begin{pf} Properties (i)-(vi) are proved by duplicating the proof of \cite[Theorem 3.2]{Cos}, so we will prove only (vii)-(ix).

\vskip 2 mm

 In order to prove (vii) and (viii), it is enough to assume that $k=2.$ A finite induction on $k$ would prove each of these claims.
 So we assume that $k=2.$ Let $u \in Lip_{0}(\Omega_1 \cup \Omega_2)$ and let $u_i=\chi_{\Omega_i} u, i=1,2.$
 We let $v_i$ be the restriction of $u$ to $\Omega_i$ for $i=1,2.$ Then $v_i \in Lip_{0}(\Omega_i)$ for $i=1,2.$ We notice that
 $u_i$ can be regarded as the extension of $v_i$ by $0$ to $\Omega_1 \cup \Omega_2$ for $i=1,2.$ We see that $u \in W(K_1 \cup K_2, \Omega_1 \cup \Omega_2)$ if and only if $v_i \in W(K_i, \Omega_i)$ for $i=1,2.$

 Suppose first that $1\le q \le p.$ Since  $\Omega_1$ and $\Omega_2$ are disjoint and $u=u_1+u_2$ with the functions $u_i$
 supported in $\Omega_i$ for $i=1,2,$ we obtain via Proposition \ref{superadd q le p}

 \begin{eqnarray*}
 ||\nabla u||_{L^{p,q}(\Omega_1 \cup \Omega_2, m_n; {\mathbf{R}}^n)}^p &\ge& ||\nabla u_1||_{L^{p,q}(\Omega_1 \cup \Omega_2, m_n; {\mathbf{R}}^n)}^p+||\nabla u_2||_{L^{p,q}(\Omega_1 \cup \Omega_2, m_n; {\mathbf{R}}^n)}^p\\
  &=& ||\nabla v_1||_{L^{p,q}(\Omega_1, m_n; {\mathbf{R}}^n)}^p+||\nabla v_2||_{L^{p,q}(\Omega_2, m_n; {\mathbf{R}}^n)}^p.
 \end{eqnarray*}
 This proves (vii).

 Suppose now that $p < q <\infty.$ Since $\Omega_1$ and $\Omega_2$ are disjoint and $u=u_1+u_2$ with the functions $u_i$ supported in
 $\Omega_i$ for $i=1,2,$ we obtain via Proposition \ref{superadd p le q}

 \begin{eqnarray*}
 ||\nabla u||_{L^{p,q}(\Omega_1 \cup \Omega_2, m_n; {\mathbf{R}}^n)}^q &\ge& ||\nabla u_1||_{L^{p,q}(\Omega_1 \cup \Omega_2, m_n; {\mathbf{R}}^n)}^q+||\nabla u_2||_{L^{p,q}(\Omega_1 \cup \Omega_2, m_n; {\mathbf{R}}^n)}^q\\
 &=& ||\nabla v_1||_{L^{p,q}(\Omega_1, m_n; {\mathbf{R}}^n)}^q+||\nabla v_2||_{L^{p,q}(\Omega_2, m_n; {\mathbf{R}}^n)}^q.
 \end{eqnarray*}
 This proves (viii).

 We see that (ix) follows from (vii) and (ii) when $1 \le q \le p.$ (We use (vii) with $k=2$ by taking $K_1=K$ and $K_2=\emptyset.$) When $p < q <\infty,$ (ix) follows from (viii) and (ii). (We use (viii) with $k=2$ by taking $K_1=K$ and $K_2=\emptyset.$) This finishes the proof of the theorem.

 \end{pf}

 \begin{Remark} \label{cap of a cpt set equal cap of its bdry} The definition of the $p,q$-capacitance implies
 $${\mathrm{cap}}_{p,q}(K, \Omega)={\mathrm{cap}}_{p,q}(\partial K, \Omega)$$
 whenever $K$ is a compact set in $\Omega.$ Moreover, if $n=1$ and $\Omega$ is an open
 interval of ${\mathbf{R}},$ then
 $${\mathrm{cap}}_{p,q}(K, \Omega)={\mathrm{cap}}_{p,q}(H, \Omega),$$
 where $H$ is the smallest compact interval containing $K.$

 \end{Remark}

 \section{\bf{Conductor inequalities}}

 \begin{Lemma} \label{upper semicont of cap(Mat Mt)} Suppose $\Omega \subset {\mathbf{R}}^n$ is open. Let $f \in Lip_{0}(\Omega)$ and let $a>1$ be a constant.
 For $t > 0$ we denote $M_{t}= \{ x \in \Omega: |f(x)|>t \}.$ Then the function $t \mapsto {\mathrm{cap}}_{p,q}(\overline{M_{at}}, M_t)$
 is upper semicontinuous.
 \end{Lemma}

 \begin{pf} Let $t_0>0$ and $\varepsilon>0.$ Let $u \in W(\overline{M_{at_{0}}}, M_{t_{0}})$ be chosen such that
 $$||\nabla u||_{L^{p,q}(\Omega, m_n; {\mathbf{R}}^n)}^p < {\mathrm{cap}}_{p,q}(\overline{M_{at_{0}}}, M_{t_0})+ \varepsilon.$$
 Let $g$ be an open neighborhood of $\overline{M_{at_{0}}}$ such that $u \ge 1$ on $g.$ Since $g$ contains the compact set $\overline{M_{at_0}},$
 there exists $\delta_1>0$ small such that $g \supset \overline{M_{a(t_0-\delta_1)}}.$ Let $G$ be an open set such that
 $\mbox{supp } u \subset  G \subset \subset M_{t_0}.$ There exists a small $\delta_2>0$ such that $\overline{G} \subset M_{t_0+\delta_2}.$
 Thus we have $\overline{M_{a(t_0-\delta)}} \subset g$ and $\overline{G} \subset {M_{t_0+\delta}}$ for every $\delta \in (0, \min\{\delta_1, \delta_2\}).$ From the choice of $g$ and $G$ we have that $u \in W(K, \Omega)$ whenever $K\subset g$ and $\overline{G} \subset \Omega.$
 This and the choice of $u$ imply that
 $${\mathrm{cap}}_{p,q}(\overline{M_{a(t_0-\delta)}}, M_{t_0+\delta}) \le {\mathrm{cap}}_{p,q}(\overline{M_{at_0}}, M_{t_0})+ \varepsilon$$
 for every $\delta \in (0, \min \{ \delta_1, \delta_2 \}).$ Using the monotonicity of ${\mathrm{cap}}_{p,q},$ we deduce that
  $$ {\mathrm{cap}}_{p,q}(\overline{M_{at}}, M_{t}) \le {\mathrm{cap}}_{p,q}(\overline{M_{at_0}}, M_{t_0})+ \varepsilon$$
 for every $t$ sufficiently close to $t_0.$ The result follows.
 \end{pf}

 \begin{Theorem} \label{strong cap ineq with Phi convex}
 Let $\Phi$ denote an increasing convex (not necessarily strictly convex) function given on $[0, \infty),$ $\Phi(0)=0.$
 Suppose $a>1$ is a constant.

 \par{\rm(i)} If $1 \le q \le p,$ then
 $$\Phi^{-1}\left( \int_{0}^{\infty} \Phi(t^p {\mathrm{cap}}_{p,q}(\overline{M_{at}}, M_t)) \frac{dt}{t}   \right) \le c(a,p,q) ||\nabla \varphi||_{L^{p,q}(\Omega, m_n; {\mathbf{R}}^n)}^p$$
 for every $\varphi \in Lip_{0}(\Omega).$

 \par{\rm(ii)} If $p < q<\infty,$ then
 $$\Phi^{-1}\left( \int_{0}^{\infty} \Phi(t^q {\mathrm{cap}}_{p,q}(\overline{M_{at}}, M_t)^{q/p}) \frac{dt}{t}   \right) \le c(a,p,q) ||\nabla \varphi||_{L^{p,q}(\Omega, m_n; {\mathbf{R}}^n)}^q$$
 for every $\varphi \in Lip_{0}(\Omega).$
 \end{Theorem}

 \begin{pf} The proof follows \cite{M4}. When $p=q$ we are in the case of the $p$-capacitance and for that case the result was proved in
 \cite[Theorem 1]{M4}. So we can assume without loss of generality that $p \neq q.$ Let $\varphi \in Lip_{0}(\Omega).$ We set
 $$\Lambda_t(\varphi)=\frac{1}{(a-1)t}  \min \{ (|\varphi|-t)_{+}, (a-1)t   \}.$$
 From Lemma \ref{same capacity if Lip and zero on bdry considered} we notice that
 \begin{equation} \label{expression for grad Lambdat}
 \Lambda_t(\varphi) \in W_{1}(\overline{M_{at}}, M_t) \mbox{ and } |\nabla \Lambda_t(\varphi)|=\frac{1}{(a-1)t} \chi_{M_t \setminus M_{at}} |\nabla \varphi| \mbox{ $m_n$-a.e.}
 \end{equation}
 The proof splits now, depending on whether $1 \le q < p$ or $p < q <\infty.$

 We assume first that $1 \le q < p.$
 From (\ref{expression for grad Lambdat}) we have
 $$t^{p} {\mathrm{cap}}_{p,q}(\overline{M_{at}}, M_t) \le \frac{1}{(a-1)^p} || \chi_{M_t \setminus M_{at}} \nabla \varphi||_{L^{p,q}(\Omega, m_n; {\mathbf{R}}^n)}^p.$$
 Hence
 \begin{equation*} \label{condcap q le p}
 \int_{0}^{\infty} \Phi(t^p {\mathrm{cap}}_{p,q}(\overline{M_{at}}, M_t)) \frac{dt}{t} \le
 \int_{0}^{\infty} \Phi(\frac{1}{(a-1)^p} ||\chi_{M_t \setminus M_{at}} \nabla \varphi||_{L^{p,q}(\Omega, m_n; {\mathbf{R}}^n)}^p) \frac{dt}{t}.
 \end{equation*}
 Let $\gamma$ denote a locally integrable function on $(0, \infty)$ such that there exist the limits $\gamma(0)$ and $\gamma(\infty).$ Then the identity
 \begin{equation} \label{integral formula involving gamma}
 \int_{0}^{\infty} (\gamma(t)-\gamma(at)) \frac{dt}{t} = (\gamma(0)-\gamma(\infty)) \log a
 \end{equation}
 holds.

 We set
 $$\gamma(t)= \Phi(\frac{1}{(a-1)^p} || \chi_{M_t} \nabla \varphi||_{L^{p,q}(\Omega, m_n; {\mathbf{R}}^n)}^p ).$$
 Using the monotonicity and convexity of $\Phi$ together with Proposition \ref{superadd q le p} and the definition of $\gamma,$ we see that
 $$\Phi( \frac{1}{(a-1)^p} || \chi_{M_t \setminus M_{at}} \nabla \varphi||_{L^{p,q}(\Omega, m_n; {\mathbf{R}}^n)}^p ) \le \gamma(t)-\gamma(at)
 \mbox{ for every $t>0$}.$$
 Since $$\gamma(0)=\Phi(\frac{1}{(a-1)^p} ||\nabla \varphi||_{L^{p,q}(\Omega, m_n; {\mathbf{R}}^n)}^p ) \mbox{ and } \gamma(\infty)=0,$$
 we get
  $$\int_{0}^{\infty} \Phi(t^p {\mathrm{cap}}_{p,q}(\overline{M_{at}}, M_t)) \frac{dt}{t} \le \log a \cdot \Phi(\frac{1}{(a-1)^p}
 ||\nabla \varphi||_{L^{p,q}(\Omega, m_n; {\mathbf{R}}^n)}^p).$$ This finishes the proof of the case $1 \le q < p.$

 We assume now that $p<q<\infty.$ From (\ref{expression for grad Lambdat}) we have
 $$t^{q} {\mathrm{cap}}_{p,q}(\overline{M_{at}}, M_t)^{q/p} \le \frac{1}{(a-1)^q} || \chi_{M_t \setminus M_{at}} \nabla \varphi||_{L^{p,q}(\Omega, m_n; {\mathbf{R}}^n)}^q.$$
  Hence
 \begin{equation*}
 \int_{0}^{\infty} \Phi(t^q {\mathrm{cap}}_{p,q}(\overline{M_{at}}, M_t)^{q/p}) \frac{dt}{t} \le
 \int_{0}^{\infty} \Phi(\frac{1}{(a-1)^q} ||\chi_{M_t \setminus M_{at}} \nabla \varphi||_{L^{p,q}(\Omega, m_n; {\mathbf{R}}^n)}^q) \frac{dt}{t}.
 \end{equation*}
 As before, we let $\gamma$ denote a locally integrable function on $(0, \infty)$ such that there exist the limits $\gamma(0)$ and $\gamma(\infty).$ We set
 $$\gamma(t)= \Phi(\frac{1}{(a-1)^q} || \chi_{M_t} \nabla \varphi||_{L^{p,q}(\Omega, m_n; {\mathbf{R}}^n)}^q ).$$
 Using the monotonicity and convexity of $\Phi$ together with Proposition \ref{superadd p le q} and the definition of $\gamma,$ we see that
 $$\Phi( \frac{1}{(a-1)^q} || \chi_{M_t \setminus M_{at}} \nabla \varphi||_{L^{p,q}(\Omega, m_n; {\mathbf{R}}^n)}^q ) \le \gamma(t)-\gamma(at)
 \mbox{ for every $t>0$}.$$
 Since $$\gamma(0)=\Phi(\frac{1}{(a-1)^q} ||\nabla \varphi||_{L^{p,q}(\Omega, m_n; {\mathbf{R}}^n)}^q ) \mbox{ and }
 \gamma(\infty)=0,$$ we get

 $$\int_{0}^{\infty} \Phi(t^q {\mathrm{cap}}_{p,q}(\overline{M_{at}}, M_t)^{q/p}) \frac{dt}{t} \le \log a \cdot \Phi(\frac{1}{(a-1)^q}
 ||\nabla \varphi||_{L^{p,q}(\Omega, m_n; {\mathbf{R}}^n)}^q).$$ This finishes the proof of the case $p < q<\infty.$
 The theorem is proved.

 \end{pf}

 Choosing $\Phi(t)=t,$ we arrive at the inequalities mentioned in the beginning of this paper.

 \begin{Corollary} \label{strong cap ineq with Phi equal identity}
 Suppose $1<p<\infty$ and $1 \le q<\infty.$ Let $a>1$ be a constant. Then (\ref{stronger cap Sobolev-Lorentz q le p}) and
 (\ref{stronger cap Sobolev-Lorentz p le q}) hold for every $\varphi \in Lip_{0}(\Omega).$
 \end{Corollary}

 \section{\bf{Necessary and sufficient conditions for two-weight embeddings}}

 We derive now necessary and sufficient conditions for Sobolev-Lorentz type inequalities involving two measures, generalizing results
 obtained in \cite{M4} and \cite{M5}.

 \begin{Theorem} Let $p,q,r,s$ be chosen such that $1<p<\infty$, $1 \le q <\infty$ and $1<s \le \max(p,q) \le r <\infty.$
 Let $\Omega$ be an open set in ${\mathbf{R}}^n$ and let $\mu$ and $\nu$ be two nonnegative locally finite measures on $\Omega.$

 \par{\rm(i)} Suppose that $1 \le q \le p.$
 The inequality
 \begin{equation} \label{two weight Lorentz form q le p}
  ||f||_{L^{r,p}(\Omega, \mu)} \le A \left( ||\nabla f||_{L^{p,q}(\Omega, m_n; {\mathbf{R}}^n)} + ||f||_{L^{s,p}(\Omega, \nu)}  \right)
 \end{equation}
 holds for every $f \in Lip_{0}(\Omega)$ if and only if there exists a constant $K>0$ such that the inequality
 (\ref{two measure inequality for general n}) is valid for all open bounded sets $g$ and $G$ that are subject to $\overline{g} \subset G \subset \overline{G} \subset \Omega.$

 \par{\rm(ii)} Suppose that $p<q<\infty.$ The inequality
 \begin{equation} \label{two weight Lorentz form p le q}
 ||f||_{L^{r,q}(\Omega, \mu)} \le A \left( ||\nabla f||_{L^{p,q}(\Omega, m_n; {\mathbf{R}}^n)} + ||f||_{L^{s,q}(\Omega, \nu)}  \right)
 \end{equation}
 holds for every $f \in Lip_{0}(\Omega)$ if and only if there exists a constant $K>0$ such that the inequality
 (\ref{two measure inequality for general n}) is valid for all open bounded sets $g$ and $G$ that are subject to $\overline{g} \subset G \subset \overline{G} \subset \Omega.$

 \end{Theorem}

 \begin{pf} We suppose first that $1 \le q \le p.$ The case $q=p$ was studied in \cite{M5}. Without loss of generality
 we can assume that $q<p.$ We choose some bounded open sets $g$ and $G$ such that
 $\overline{g} \subset G \subset \overline{G} \subset \Omega$ and $f \in W(\overline{g}, G)$ with $0 \le f \le 1.$
 We have
 \begin{equation*}
 \mu(g) \le C(r,p) \, ||f||_{L^{r,p}(\Omega, \mu)}^r
 \end{equation*}
 and
 \begin{equation*}
 ||f||_{L^{s,p}(\Omega, \nu)}^s \le C(s,p) \, \nu(G)
 \end{equation*}
 for every $f \in W(\overline{g}, G)$ with $0 \le f \le 1.$ The necessity for $1 \le q<p$ is obtained by taking the infimum over all such functions $f$ that are admissible for the conductor $(\overline{g}, G).$

 We prove the sufficiency now when $1 \le q<p$. Let $a \in (1, \infty).$
 We have
 $$a^p \int_{0}^{\infty} \mu(M_{at})^{p/r} d(t^p) \le a^p K_1\left(\int_{0}^{\infty} ({\mathrm{cap}}_{p,q}(\overline{M_{at}}, M_{t}) + \nu(M_t)^{p/s}) d(t^p) \right).$$
 This and (\ref{stronger cap Sobolev-Lorentz q le p}) yield the sufficiency for the case $1 \le q < p.$

 Suppose now that $p<q<\infty.$ We choose some bounded open sets $g$ and $G$ such that
 $\overline{g} \subset G \subset \overline{G} \subset \Omega$ and $f \in W(\overline{g}, G)$ with $0 \le f \le 1.$
 We have
 \begin{equation*}
 \mu(g) \le C(r,q) \, ||f||_{L^{r,q}(\Omega, \mu)}^r
 \end{equation*}
 and
 \begin{equation*}
 ||f||_{L^{s,q}(\Omega, \nu)}^s \le C(s,q) \, \nu(G)
 \end{equation*}
 for every $f \in W(\overline{g}, G)$ with $0 \le f \le 1.$ The necessity for $p < q<\infty$ is obtained by taking the infimum over all such functions $f$ that are admissible for the conductor $(\overline{g}, G).$

 We prove the sufficiency now when $p < q<\infty.$ Let $a \in (1, \infty).$
 We have
 $$a^q \int_{0}^{\infty} \mu(M_{at})^{q/r} d(t^q) \le a^q K_2 \left(\int_{0}^{\infty} ({\mathrm{cap}}_{p,q}(\overline{M_{at}}, M_t)^{q/p} + \nu(M_t)^{q/s}) d(t^q) \right).$$
 This and (\ref{stronger cap Sobolev-Lorentz p le q}) yield the sufficiency for the case $p < q<\infty.$ The proof is finished.
 \end{pf}

 We look for a simplified necessary and sufficient two-weight imbedding condition when $n=1.$ Before we state and prove such a condition
 for the case $n=1,$ we need to obtain sharp estimates for the $p,q$-capacitance of conductors $([a,b], (A,B))$ with $A<a<b<B.$ This is
 the goal of the following proposition.

 \begin{Proposition} \label{sharp estimates for capacity of intervals when n is 1}
 Suppose $n=1,$ $1<p<\infty$ and $1 \le q \le \infty.$ There exists a constant $C(p,q)\ge 1$ such that
 $$C(p,q)^{-1} (\sigma_{1}^{1-p}+ \sigma_{2}^{1-p}) \le {\mathrm{cap}}_{p,q}([a,b], (A, B)) \le C(p,q) (\sigma_{1}^{1-p}+ \sigma_{2}^{1-p}),$$
 where $\sigma_1=a-A$ and $\sigma_2=B-b.$

 \end{Proposition}

 \begin{pf} By the behaviour of the Lorentz $p,q$-quasinorm in $q$ (see for instance \cite[Proposition IV.4.2]{BS}), it suffices to find
 the upper bound for the $p,1$-capacitance and the lower bound for the $p, \infty$-capacitance of the conductor $([a,b], (A,B)).$
 We start with the upper bound for the $p,1$-capacitance of this conductor.

 We use the function $u: (A, B) \rightarrow {\mathbf{R}}$ defined by

 $$u(x)=\left\{ \begin{array}{cl}
  1 & \mbox{if $a \le x \le b$} \\
 \frac{x-A}{\sigma_1} & \mbox{if $A<x<a$}\\
 \frac{B-x}{\sigma_2} & \mbox{if $b<x<B.$}
 \end{array}
 \right.$$

 Then from Lemma \ref{same capacity if Lip and zero on bdry considered} it follows that $u \in W_{1}([a,b], (A, B))$
 with
 $$|u'(x)|=\left\{ \begin{array}{cl}
  0 & \mbox{if $a < x < b$} \\
 {\sigma_1}^{-1} & \mbox{if $A<x<a$}\\
 {\sigma_2}^{-1} & \mbox{if $b<x<B.$}
 \end{array}
 \right.$$

 We want to compute an upper estimate for $||u'||_{L^{p,1}((A,B), m_1)}.$ We have
 \begin{eqnarray}
 ||u'||_{L^{p,1}((A,B), m_1)} &\le& ||{\sigma_1}^{-1}||_{L^{p,1}((A,a), m_1)} +
 ||{\sigma_2}^{-1}||_{L^{p,1}((b,B), m_1)} \\
 &=& p \left( \sigma_1^{-1+1/p}+ \sigma_2^{-1+1/p} \right) \nonumber.
 \end{eqnarray}

 Therefore $${\mathrm{cap}}_{p,1}([a,b], (A, B)) \le C(p) (\sigma_1^{1-p} + \sigma_2^{1-p}).$$

 We try to get lower estimates for the $p,\infty$-capacitance of this conductor. Let $v \in W([a,b], (A,B))$ be an arbitrary
 admissible function. We let $v_1$ be the restriction of $v$ to $(A, a)$ and $v_2$ be the restriction of $v$ to $(b, B)$
 respectively. We notice that $v'$ is supported in $(A, a) \cup (b, B).$ Therefore, since $v'$ coincides with
 $v_1'$ on $(A, a)$ and with $v_2'$ on $(b, B),$ we have that
 \begin{equation} \label{lower bound for p infty capacity when n is one}
 ||v'||_{L^{p, \infty}((A,B), m_1)} \ge \max( ||v_1'||_{L^{p, \infty}((A, a), m_1)}, ||v_2'||_{L^{p, \infty}((b, B), m_1)}).
 \end{equation}

 From (\cite[Corollary 2.4]{Cos}) we have
 $$||v_1'||_{L^{p, \infty}((A, a), m_1)} \ge 1/p' \cdot \sigma_{1}^{-1/p'} ||v_1'||_{L^{1}((A, a), m_1)}$$
 and
 $$||v_2'||_{L^{p, \infty}((b, B), m_1)} \ge 1/p' \cdot \sigma_{2}^{-1/p'} ||v_2'||_{L^{1}((b, B), m_1)}.$$
 Since $$||v_1'||_{L^{1}((A, a), m_1)}=\int_{A}^{a} |v_1'(x)| dx \ge 1,$$
 we obtain
 \begin{equation} \label{lower bound for v1 prime p, infty}
 ||v_1'||_{L^{p, \infty}((A, a), m_1)}\ge 1/p' \cdot \sigma_1^{-1/p'}.
 \end{equation}
 Similarly, since $$||v_2'||_{L^{1}((b, B), m_1)}=\int_{b}^{B} |v_2'(x)| dx \ge 1,$$ we obtain
 \begin{equation} \label{lower bound for v2 prime p, infty}
 ||v_2'||_{L^{p, \infty}((b, B), m_1)}\ge 1/p' \cdot \sigma_2^{-1/p'}.
 \end{equation}
 From (\ref{lower bound for p infty capacity when n is one}), (\ref{lower bound for v1 prime p, infty}) and (\ref{lower bound for v2 prime p, infty}) we get the desired lower bound for the $p, \infty$-capacitance. This finishes the proof.

 \end{pf}

 Now we state and prove a necessary and sufficient two-weight imbedding condition for the case $n=1.$

 \begin{Theorem} Suppose $n=1.$ Let $p,q,r,s$ be chosen such that $1<p<\infty,$ $1 \le q < \infty$ and $1< s \le \max (p,q) \le r <\infty.$
 Let $\Omega$ be an open set in ${\mathbf{R}}$ and let $\mu$ and $\nu$ be two nonnegative locally finite measures on $\Omega.$
 \par{\rm(i)} Suppose that $1 \le q \le p.$ The inequality
 \begin{equation} \label{two weight Lorentz form q le p n is 1}
  ||f||_{L^{r,p}(\Omega, \mu)} \le A \left( ||f'||_{L^{p,q}(\Omega, m_1)} + ||f||_{L^{s,p}(\Omega, \nu)}  \right)
 \end{equation}
 holds for every $f \in Lip_{0}(\Omega)$ if and only if there exists a constant $K>0$ such that the inequality
 (\ref{two measure inequality for centered intervals when n equal 1}) is valid whenever $x,$ $d$ and $\tau$ are such that
 $\overline{\sigma_{d+\tau}(x)} \subset \Omega.$

 \par{\rm(ii)} Suppose that $p < q <\infty.$ The inequality
 \begin{equation} \label{two weight Lorentz form p le q n is 1}
  ||f||_{L^{r,q}(\Omega, \mu)} \le A \left( ||f'||_{L^{p,q}(\Omega, m_1)} + ||f||_{L^{s,q}(\Omega, \nu)}  \right)
 \end{equation}
 holds for every $f \in Lip_{0}(\Omega)$ if and only if there exists a constant $K>0$ such that the inequality
 (\ref{two measure inequality for centered intervals when n equal 1}) is valid whenever $x,$ $d$ and $\tau$ are such that
 $\overline{\sigma_{d+\tau}(x)} \subset \Omega.$
  \end{Theorem}

 \begin{pf} We only have to prove that the sufficiency condition for intervals implies the sufficiency condition for general bounded and open sets
 $g$ and $G$ with $\overline{g} \subset G \subset \overline{G} \subset \Omega.$ Let $G$ be the union of nonoverlapping intervals $G_i$
 and let $g_i=G \cap g_i.$ We denote by $h_i$ the smallest interval containing $g_i$ and by $\tau_{i}$ the minimal distance from $h_i$ to ${\mathbf{R}} \setminus G_i.$ We also denote by $H_i$ the open interval concentric with $h_i$ such that the minimal distance from $h_i$ to ${\mathbf{R}} \setminus H_i$ is $\tau_i.$ Then $H_i \subset G_i.$
 From Remark \ref{cap of a cpt set equal cap of its bdry} we have that ${\mathrm{cap}}_{p,q}(\overline{g_i}, G_i)={\mathrm{cap}}_{p,q}(\overline{h_i}, G_i).$ Moreover, from Theorem \ref{Cap Thm 1<q leq infty} (ii) and Proposition \ref{sharp estimates for capacity of intervals when n is 1} we have
 $$C(p,q)^{-1} \tau_{i}^{1-p} \le {\mathrm{cap}}_{p,q}(\overline{h_i}, G_i) \le {\mathrm{cap}}_{p,q}(\overline{h_i}, H_i) \le 2 \, C(p,q) \tau_{i}^{1-p}$$ for some constant $C(p,q) \ge 1.$
 Since $\overline{g}$ is compact lying in $\cup_{i \ge 1} G_i,$ it follows that $\overline{g}$ is covered by
 only finitely many of the sets $G_i.$ This and Theorem \ref{Cap Thm 1<q leq infty} (ix) allow us to assume
 that $G$ is in fact written as a finite union of disjoint intervals $G_i.$ Now the proof splits, depending
 on whether $1 \le q \le p$ or $p < q<\infty.$

 We assume first that $1 \le q \le p.$
 We have
 \begin{equation} \label{superadd for p,q cap q le p}
 {\mathrm{cap}}_{p,q}(\overline{g}, G) \ge \sum_{i} {\mathrm{cap}}_{p,q}(\overline{g_i}, G_i)= \sum_{i} {\mathrm{cap}}_{p,q}(\overline{h_i}, G_i).
 \end{equation}

 Using (\ref{two measure inequality for centered intervals when n equal 1}), we obtain
 \begin{eqnarray*}
 \mu(g_i)^{p/r} &\le& \mu(h_i)^{p/r} \le K_1 (\tau_i^{1-p}+ \nu(H_i)^{p/s}) \\
 &\le& K_1 \, C(p,q) ({\mathrm{cap}}_{p,q}(\overline{g_i}, G_i)+ \nu(G_i)^{p/s})
 \end{eqnarray*}
 where $K_1$ is a positive constant independent of $g$ and $G.$ Since $s \le p \le r<\infty,$ we have
 $$\mu(g)^{p/r} \le \sum_{i} \mu(g_i)^{p/r}$$
 and
 $$\sum_{i} \nu(G_i)^{p/s} \le \nu(G)^{p/s}.$$
 This and (\ref{superadd for p,q cap q le p}) prove the claim when $1 \le q \le p.$

 We assume now that $p<q<\infty.$ We have
 \begin{equation} \label{superadd for p,q cap p le q}
 {\mathrm{cap}}_{p,q}(\overline{g}, G)^{q/p} \ge \sum_{i} {\mathrm{cap}}_{p,q}(\overline{g_i}, G_i)^{q/p}= \sum_{i} {\mathrm{cap}}_{p,q}(\overline{h_i}, G_i)^{q/p}.
 \end{equation}

 Using (\ref{two measure inequality for centered intervals when n equal 1}), we obtain
 \begin{eqnarray*}
 \mu(g_i)^{q/r} &\le& \mu(h_i)^{q/r} \le K_2 (\tau_i^{q(1-p)/p}+ \nu(H_i)^{q/s}) \\
 &\le& K_2 \, C(p,q)^{q/p} ({\mathrm{cap}}_{p,q}(\overline{g_i}, G_i)^{q/p}+ \nu(G_i)^{q/s})
 \end{eqnarray*}
 where $K_2$ is a positive constant independent of $g$ and $G.$ Since $s \le q \le r<\infty,$ we have
 $$\mu(g)^{q/r} \le \sum_{i} \mu(g_i)^{q/r}$$
 and $$\sum_{i} \nu(G_i)^{q/s} \le \nu(G)^{q/s}.$$
 This and (\ref{superadd for p,q cap p le q}) prove the claim when $p<q<\infty.$ The theorem is proved.
\end{pf}
{\bf{Acknowledgements.}} Both authors were partially supported by NSERC and by the Fields Institute (\c S. Costea) and by NSF grant DMS 0500029
(V. Maz'ya).

\end{document}